\documentclass[12pt]{article}
\usepackage{amsmath,amsfonts}
\begin{document}

\baselineskip 20pt
\date{}
\title{On the theorem of F. Riesz on Blaschke products and on inner functions in the unit ball of 
$\mathbb C^{n}$\footnote{2010 Mathematics Subject Classification: Primary 30J05, 30J10, 32A40. Key words:
Blaschke products, inner functions in the unit ball of $\mathbb C^{n}$, $G_\delta$ sets, radial limit.}}

\author{Arthur A.~Danielyan}

\maketitle

\begin{abstract}
 
\noindent For a convergent Blaschke product $f$ in the unit disc of $\mathbb C$,
let $G(f)$ be the set of all points $t$ on the unit circle for which 
$|f(rt)| \rightarrow 1$ as $r \rightarrow 1$.  
By a classical F. Riesz theorem 
$G(f)$ is of full measure, which describes $G$ metrically. 
We describe $G$ topologically by proving that it is a $G_{\delta}$ set. 
The same description is true in general for inner functions not only in the unit disc but also in the unit ball of $\mathbb C^n$.

\end{abstract}

\begin{section}{Introduction.}

Let  $U = \{z: |z|<1\}$ and $T= \{z: |z|=1\}$ be the open unit disc and the unit circle in the complex plane $\mathbb C$, respectively.
We denote by $m$ the Lebesgue measure on $T$. 
 We also denote by $B$ and by $S$ the unit ball and the unit sphere of $\mathbb C^n$, respectively. 
 (When $n=1$, then $B=U$ and $S=T$; cf. \cite{Rud}, p. 2.) 
 On $S$ the standard measure is the unique rotation-invariant positive Borel measure $\sigma$ for which $\sigma(S)=1$ 
 (see \cite{Rud}, p. 12, or \cite{Rud2}, p. 100). 
 
 For a sequence $\{a_k\}$ of points in $U$ with $ \sum_{k=1}^{\infty} (1 - |a_k|) < \infty$, the  Blaschke product
$$f(z)=\prod_{k=1}^{\infty}\frac{|a_k|}{a_k}\frac{a_k-z}{1 -\bar a_k z} $$
is a bounded (by 1) analytic function in $U$ the zeros of which are exactly at the points $a_k$.

For a function $f$ defined in $U$  (respectively in $B$) and for $\zeta \in T$ (respectively $\zeta \in S$), we denote
$f(\zeta) =  \lim_{r \to 1} f(r\zeta),$ if the limit exists.

By Fatou's theorem, the radial limits of any bounded analytic function exist a.e. on $T$. 
Moreover, the following is true.    

  \vspace{0.4 cm}

{\bf Theorem A (F. Riesz, 1923).}\ {\it The radial limits of a Blaschke product are of modulus 1 a.e. on $T$.}
 
  \vspace{0.4 cm}
 
The inner functions are making a subclass of bounded analytic functions. One can define them both for $\mathbb C$
and $\mathbb C^n$ ($n \geq 2$) as follows (cf. \cite{Hoff}, p. 62, and  \cite{Rud2}, p. 101, \cite{Rud}, p. 403). 

 \vspace{0.4 cm}
 
 {\bf Definition 1.}  { \it A nonconstant function $f \in H^{\infty}(B)$, $n \geq 1$, is an inner function if the radial limits, $\lim_{r \rightarrow 1}  f(rt) = f(t),$ 
 satisfy $|f(t)|=1$ for almost all $t \in S.$}
 
  \vspace{0.4 cm}
  
  By Theorem A, the Blaschke products are inner functions defined in $U$. However, it is well known that  
  there are other inner functions too in $U$. 
  The existence of inner functions in $B$ of $\mathbb C^n$  has been proved by Aleksandrov \cite{Alex} (cf. \cite{Rud2}). 
 
  Obviously in Definition 1 the condition $f \in H^{\infty}(B)$ can be replaced by the condition $|f(z)| \leq 1,  z \in B$ (cf. \cite{Alex}, p. 143).
  But since $f$ is nonconstant, the maximum modulus theorem (cf. \cite{Rud}, p. 5, (iv)) implies that in fact we have  $|f(z)| <1,  z \in B$. 
  
  \vspace{0.4 cm}  
  
{\bf Definition 2.} {\it For an inner function $f$ in $B$ of $\mathbb C^n$, $n \geq1$, define the set

$$G(f)=\{t \in S:  \lim_{r \to 1} |f(rt)|=1\}.$$}

  \vspace{0.4 cm}

For a Blaschke product $f$ in $U$, by Theorem A, $G(f)$
is of full measure, which is a metric description of $G(f)$.
 Obviously, $G(f)$ is of full measure for any inner function in $B$ of $\mathbb C^{n}$, for each $n \geq 1.$

The purpose of the present paper is to describe the set $G(f)$ topologically. 

  \vspace{0.4 cm}
 
{\bf Theorem 1.}\ {\it If $f$ is an inner function in $B$ of $\mathbb C^n$, $n \geq 1,$
then  $G(f)$ is a $G_\delta$ set.}
 
  \vspace{0.4 cm}
  
  When $n=1$, Theorem 1 is a topological counterpart of Theorem A. 
  We formulate this explicitly for Blaschke products as a corollary.

  \vspace{0.4 cm}

{\bf Corollary 1.}\ {\it If $f$ is a Blaschke product in $U$, then $G(f)$ is a $G_\delta$ set.}
 
  \vspace{0.4 cm}

 The following corollaries of Theorem 1 are also obvious.  
 
  \vspace{0.4 cm} 
   
{\bf Corollary 2.}\ {\it If $f$ is an inner function in $U$, then $T \setminus G(f)$ is $F_\sigma$ and $m(T \setminus G(f))=0.$}
 
  \vspace{0.4 cm}

  \vspace{0.4 cm}
 
{\bf Corollary 3.}\ {\it If $f$ is an inner function in $B$, then $S \setminus G(f)$ is $F_\sigma$ and $\sigma(S \setminus G(f))=0$.}
 
  \vspace{0.4 cm}

Since an $F_\sigma$ set of measure zero is meagre\footnote{A meagre set is a countable union of nowhere dense sets (see \cite{Car}).},
Corollary 2 implies the following known results of G. T. Cargo \cite{Car} (see Theorem 1 and Theorem 2
in  \cite{Car}) which we  formulate as Theorem B and Theorem C.

   \vspace{0.4 cm}
 
{\bf Theorem B.}\ {\it If $f$ is an inner function in $U$, then the set   $\{\zeta \in T :  f(\zeta)=\alpha\}$ is 
meagre in $T$ for each $\alpha \in U$.}
 
  \vspace{0.4 cm} 
  
  {\bf Theorem C.}\ {\it If $f$ is an inner function in $U$, then  $\{\zeta \in T :  f(\zeta) \in U \}$ is meagre in $T$.}  
  
  \vspace{0.4 cm}

 The original proofs of the formulated theorems are based on some well-known results due to Seidel, Bagemihl and Seidel, and Collingwood 
 from the theory of the cluster sets and related areas (see \cite{Car}).  
 In contrast, the proof of the present paper is elementary and straightforward.

Corollary 3 extends Theorem B and Theorem C to several complex variable setting. 
As far as we know, not only Corollary 3, but also the mentioned extensions are new results in several complex variables.  
Also, the present paper suggests a unified method of proof for such results both 
for one and several complex variables cases.  
Since Theorem C implies Theorem B, we only formulate the extension of Theorem C to the case of several complex variables. 
 \vspace{0.4 cm} 
  
  {\bf Corollary 4.}\ {\it If $f$ is an inner function in $B$, then  $\{\zeta \in S :  f(\zeta) \in U \}$ is meagre in $S$.}  
  
  \vspace{0.4 cm}

  \vspace{0.4 cm}

\end{section}

\begin{section}{Proofs.}

The proofs of the corollaries are obvious. We present the proof of Theorem 1.

 {\it Proof of Theorem 1}.  Let $f$ be  an inner function in $B$ of $\mathbb C^n$, $n \geq 1.$ 
For $t \in S$, obviously $ \lim_{r \to 1} |f(rt)|=1$ if and only if $ \liminf_{r \to 1} |f(rt)|=1.$ This implies that $t \in G(f)$ if and only if 
$$ \limsup_{r \to 1} \frac{1 }{1-|f(rt)|}=+\infty.$$
In other words,  
$$G(f)= \{t\in S:  \limsup_{r \to 1} \frac{1 }{1-|f(rt)|}=+\infty \}.$$
Such a set is $G_\delta$. We present a direct and short proof for it. 

The nonnegative function 
$$ \varphi_r(t)= \frac{1 }{1-|f(rt)|}$$ 
is continuous for fixed $r$, $0<r<1$, as a function of $t$ on $S$.  
Denote $D^{(k)}_{r} =\{t:  \varphi_r(t) \leq k\}.$ Since $D^{(k)}_{r}$ is closed for each $r$, the intersection 
$$D^{(k)}= \bigcap_{0<r<1} D^{(k)}_{r}$$
 is closed.  For $t_0 \in D^{(k)}$,  $ \varphi_r(t_0) \leq k$ for all $r$. Thus 
$ \limsup_{r \to 1} \varphi_r(t_0) \leq k,$ which means that $t_0 \in S \setminus G(f)$. Thus
$$\bigcup_{k=1}^{\infty} D^{(k)} \subset S \setminus G(f). $$

Also 
$$ S \setminus G(f)  \subset \bigcup_{k=1}^{\infty} D^{(k)},$$
as if $ \limsup_{r \to 1} \varphi_r(t_0) < +\infty,$ there exists $k$ such that $ \varphi_r(t_0) \leq k$ for all $0<r<1$, and so $t_0 \in D^{(k)}$. 
Hence,
$$ S \setminus G(f) = \bigcup_{k=1}^{\infty} D^{(k)}.$$

Thus $ S \setminus G(f)$ is $F_\sigma$, which means $G(f)$ is $G_\delta$.  Theorem 1 is proved. 
 \end{section}

 \vspace{1 cm}
 
{\bf AI Declaration.} AI was not used in preparation of this manuscript.

\begin{minipage}[t]{6.5cm}
Arthur A. Danielyan\\
Department of Mathematics\\ 
and Statistics\\
University of South Florida\\
Tampa, Florida 33620\\
USA\\
{\small e-mail: adaniely@usf.edu}\\
\end{minipage}

\end{document}